\documentclass[12pt]{article}
\usepackage{amsmath}
\usepackage{amssymb}
\usepackage{graphicx}
\usepackage{grffile}
\usepackage{enumerate}
\usepackage{natbib}
\usepackage{url} 

\usepackage[usenames,dvipsnames]{xcolor}

\newcommand{\blind}{0}

\newtheorem{Assum}{Assumption}
\newtheorem{Theor}{Theorem}
\newtheorem{Corol}{Corollary}[section]

\newcommand{\cqfd}{\hfill $\square$}

\newcommand{\R}{\mathbb R}

\begin{document}

\def\spacingset#1{\renewcommand{\baselinestretch}%
{#1}\small\normalsize} \spacingset{1}


\if0\blind
{
  \title{\bf Revisiting the name variant of\\ the two-children problem}
  \author{Davy Paindaveine 
\\
    ECARES and Department of Mathematics, Universit\'{e} libre de Bruxelles\\
    and \\
    Philippe Spindel \\
    Service de Physique de l'Univers, Champs et Gravitation, Universit\'{e} de Mons\\
Service de Physique Th\'{e}orique, Universit\'{e} libre de Bruxelles}
  \maketitle
} \fi

\if1\blind
{
  \bigskip
  \bigskip
  \bigskip
  \begin{center}
    {\LARGE\bf Revisiting the name variant of\\ the two-children problem}
\end{center}
  \medskip
} \fi

\bigskip
\begin{abstract}
Initially proposed by Martin Gardner in the 1950s, the famous two-children problem is often presented as a paradox in probability theory. A relatively recent variant of this paradox states that, while in a two-children family for which at least one child is a girl, the probability that the other child is a boy is~$2/3$, this probability becomes~$1/2$ if the first name of the girl is disclosed (provided that two sisters may not be given the same first name). We revisit this variant of the problem and show that, if one adopts a natural model for the way first names are given to girls, then the probability that the other child is a boy may take \emph{any} value in~$\textcolor{black}{(}0,2/3\textcolor{black}{)}$. By exploiting the concept of Schur-concavity, we study how this probability depends on model parameters.
\end{abstract}

\noindent%
\vfill

\newpage
\spacingset{1.45} 


\section{Introduction}
\label{sec:intro}

In one of his famous \emph{mathematical games}, Martin Gardner asked the following questions: \emph{``\mbox{Mr.} Smith has two children. At least one of them is a boy. What is the probability that the other child is a boy? \mbox{Mr.} Jones has two children. The older is a girl. What is the probability that the other child is a girl?"} See \cite{Gardner1959}. While the answers he first provided were~$1/3$ and~$1/2$, respectively, he explicitly stated later that the first question was actually ambiguous; see \cite{Gardner1987}, Chapters~14 and~19. To phrase this question without ambiguity in such a way that the corresponding answer is indeed~$1/3$, one should adopt the view (as we do throughout the present note) that \mbox{Mr.} Smith is randomly selected among all two-children families having at least one boy; see \cite{Barhillel1982} or \cite{Khovanova2011} for discussions on how the answer to this first question depends on the way one obtains the information that \mbox{Mr.} Smith has at least one boy.

The two-children problem, that is sometimes referred to as the boy-or-girl problem, gained much popularity twenty-five years ago when it was discussed by the well-known columnist Marilyn vos Savant in \emph{Parade} magazine (\citealp{vosSavant1997}). It has since been discussed in several monographs (\citealp{Mlodinow2008,Chang2012}) and scientific papers (among which \citealp{Dagostini2010,Lynch2011,Pollak2013}, and the aforementioned ones). 

Two variants of the two-children problem are famous. The first one asks: \emph{for a two-children family having at least a girl who is born on a Tuesday, what is the probability that the other child is a boy?} See \cite{Lynch2011}, \cite{Falk2011}, \cite{Taylor2014},  \cite{Zaskis2015}. If the probability that a girl is born on a Tuesday is~$r$, then this probability can be shown to be equal to~$2/(4-r)$, which ranges from~$1/2$ (for~$r=0$) to~$2/3$ (for~$r=1$); for~$r=1/7$, the probability that the other child is a boy is thus~$14/27$. This still assumes that the two-children family considered here is randomly selected from all two-children families meeting this particular property; when removing the ambiguity above in another way, different probabilities are obtained in this variant, too; see \cite{Khovanova2011}.   

The second variant rather asks: \emph{for a two-children family having at least a girl whose name\footnote{In the rest of the paper, ``name" will throughout stand for ``first name".} is Florida, what is the probability that the other child is a boy?} See, e.g., \cite{Mlodinow2008} or \cite{Marks2011}. If two sisters may be given the same name, then this variant is strictly equivalent to the previous one: \textcolor{black}{more precisely, if it is assumed that girls are independently named Florida with probability~$r$, then} the probability that the other child is a boy is~$2/(4-r)$.
To make the second variant of interest, one therefore needs to assume that two sisters may not be given the same name, in which case\textcolor{black}{, under the assumptions associated with what we will call Model~A below,} the probability that the other child is a boy is~$1/2$, irrespective of~$r$; see \cite{Dagostini2010}. 

In this note, we revisit this second variant. In Section~\ref{sec:models}, we describe two models that can be considered to answer the question of interest. The first one, Model~A, is the traditional one, for which both genders are equally likely for the other child---we carefully state the corresponding assumptions. Then, we introduce \textcolor{black}{a new} model, Model~B, that specifies the way girl names are picked by parents. In Section~\ref{sec:results}, we compute in each model the probability that the other child is a boy. While both genders are indeed equally likely in Model~A, this probability in Model~B may assume any value in~$\textcolor{black}{(}0,2/3\textcolor{black}{)}$ depending on popularities of girl names. In this second model, this probability is actually a Schur-concave function of name popularities, which allows us to study how this probability depends on model parameters. In Section~\ref{sec:finalcomments}, we provide some final comments. Finally, an appendix 
collects technical proofs.

\section{Two models}
\label{sec:models}

We consider two probabilistic models, labelled Model~A and Model~B below, for the gender and name of each child in a two-children family (in line with the question raised in the second variant of the two-children problem, we will actually consider names for girls only). We start by describing assumptions that are common to both models. Regarding gender, we adopt the following assumption.

\begin{Assum}[Models A--B]
\label{Assum1}
	(i) Any born child is a boy ($b$) with probability~$1/2$ and a girl ($g$) with probability~$1/2$.
	%
	(ii) Genders of both children are independent. 
\end{Assum}

\textcolor{black}{Writing \emph{E} and \emph{Y} for \emph{Elder} and \emph{Younger}, respectively, we thus have, with} obvious notation,~$P[Eb]=P[Eg]=P[Yb]=P[Yg]=1/2$, hence, e.g.,~$P[Eb\cap Yb]=1/4$. Now, we turn to assumptions involving names. We let~$r_1:=P[Eg_1|Eg]$, where~$Eg_1$ is the event that the elder child is a girl named~$n_1$. This event has thus probability~$P[Eg_1]=P[Eg_1|Eg]P[Eg]=r_1/2$, which yields~$P[Eg\setminus Eg_1]=1-P[Eb]-P[Eg_1]=(1-r_1)/2$. The assumptions common to both models and related to names are then as follows.

\begin{Assum}[Models A--B]
\label{Assum2}
	(i) Two sisters may not be given the same name: $P[Eg_1\cap Yg_1]=0$. 
	(ii) The name given to an elder girl and gender of the second child are independent: $P[Yb|Eg_1]=P[Yb|Eg](=1/2)$.
	(iii) A girl having an elder brother is given name~$n_1$ with the same probability as an elder girl child: $P[Yg_1|(Yg\cap Eb)]=P[Eg_1|Eg](=r_1)$.  
\end{Assum}

Assumption~\ref{Assum2}(ii) yields $P[Eg_1\cap Yb]=P[Yb|Eg_1]P[Eg_1]=r_1/4$, which implies both $P[(Eg\setminus Eg_1) \cap Yb]=P[Yb]-P[Eb\cap Yb]-P[Eg_1\cap Yb]=(1-r_1)/4$ and \textcolor{black}{(using also Assumption~\ref{Assum2}(i):)} $P[Eg_1\cap (Yg\setminus Yg_1)]=P[Eg_1]-P[Eg_1\cap Yb]\textcolor{black}{-P[Eg_1\cap Yg_1]}=r_1/4$.  
Moreover, Assumption~\ref{Assum2}(iii) provides~$P[Eb\cap Yg_1]
=P[Yg_1|(Yg \cap Eb)] P[Yg \cap Eb]
=r_1/4$, which yields~$P[Eb\cap (Yg\setminus Yg_1)]=P[Eb]-P[Eb\cap Yb]-P[Eb\cap Yg_1]=(1-r_1)/4$. 
Summing up, Assumptions~\ref{Assum1}--\ref{Assum2} lead to the probabilities given in Table~\ref{Table1}.  
Clearly, one needs an extra assumption to determine the missing probabilities in this table, and \textcolor{black}{the two models will differ}. The usual model relies on the following, often tacit, assumption.

	\begin{table}[h!]
\begin{center}
\begin{tabular}{|c|ccc|c|}
\hline
   & $Yb$ & $Yg_1$ & $Yg\setminus Yg_1$ &  \\
\hline
   $Eb$     & $1/4$ & $r_1/4$ & $(1-r_1)/4$ & $1/2$ \\
   $Eg_1$   & $r_1/4$ & $0$ & $r_1/4$ & $r_1/2$ \\
   $Eg\setminus Eg_1$ & $(1-r_1)/4$ & $p_{32}$ & $p_{33}$ & $(1-r_1)/2$ \\
\hline
    & $1/2$ & $p_{\cdot 2}$ & $p_{\cdot 3}$ & 1 \\
\hline
\end{tabular}
\caption{Probabilities obtained from Assumptions~\ref{Assum1}--\ref{Assum2}.}
\label{Table1}
\end{center}
\end{table}

\begin{Assum}[Model~A]
\label{Assum3}
	$P[Yg_1|Yg]=P[Eg_1|Eg](=r_1)$, or equivalently
	%
	%
	$P[Yg_1]=P[Eg_1](=r_1/2)$. 
\end{Assum}

From a statistical point of view, this modelling assumption essentially translates the expectation that, within two-children families, there should be roughly as many girls baring the name~$n_1$ among the younger girl children as among elder girl ones. Under this assumption, \textcolor{black}{we indeed have} $p_{\cdot 2}^A=P[Yg_1]=P[Eg_1]=r_1/2$, which \textcolor{black}{allows} us to obtain~$p_{32}^A=r_1/4$, $p_{33}^A=(1-2r_1)/4$ and~$p_{\cdot 3}^A=(1-r_1)/2$. Note that this imposes that~$r_1\leq 1/2$, a restriction we will not have in the alternative model we now describe.
 
Unlike the model above, Model~B does not rely on a statistical view but rather provides a probabilistic scheme describing how girl names are \emph{picked by parents} according to popularity. Assume that there are~$K$ names, $n_1,\ldots,n_K$ say, with respective popularity~$r_1,\ldots,r_K$, where the~$r_k$'s are positive numbers that sum up to one. If the elder child is a girl, then it will be accordingly given name~$n_k$ with probability~$r_k$. If the first child is a boy and the second one is a girl, then the $K$ names are available for this girl, which will similarly be given name~$n_k$ with probability~$r_k$. However, if the first child is a girl, named~$n_k$ say, then in case a second girl is born, this name is not available anymore (Assumption~\ref{Assum2}(i)), and parents will then naturally give this girl a name according to the rescaled probabilities associated with~$r_1,\ldots,r_{k-1},r_{k+1},\ldots,r_K$. This is formalized in the following assumption.

\setcounter{Assum}{2}
\begin{Assum}[Model~B]
Girl names available are~$n_1,\ldots,n_K$. The first girl born in a family (if any) is given name~$n_k$ with probability~$r_k$; here, $r_k\geq 0$ for any~$k=1,\ldots,K$ and~$\sum_{k=1}^K r_k=1$. If the elder girl was given name~$n_k$, then the possible second girl is given name~$n_\ell$, with probability~$r_\ell(1-\delta_{k\ell})/(1-r_k)$, where~$\delta_{k\ell}$ takes value one if~$k=\ell$ and value zero otherwise: $P[Yg_\ell|(Eg_k\cap Yg)]=r_\ell(1-\delta_{k\ell})/(1-r_k)$, $k,\ell=1,\ldots,K$, where~$Eg_k$ (resp., $Yg_k$) denotes the event that the elder (resp., younger) child is a girl named~$n_k$.
\end{Assum}

In this model where~$K$ names are available, note that Assumption~\ref{Assum2}(ii) implies that $P[Yb|Eg_k]=P[Yb|Eg](=1/2)$ for any~$k=1,\ldots,K$. Consequently, Model~$B$ yields
\begin{eqnarray*}
\lefteqn{
p_{32}^B
=
\sum_{k=2}^K
P[Eg_k\cap Yg_1]
}
\\[2mm]
& & 
\hspace{10mm}
=
\sum_{k=2}^K
P[Yg_1|(Eg_k\cap Yg)]P[Yg|Eg_k]P[Eg_k]
=
\frac{1}{4}
\sum_{k=2}^K
\frac{r_1r_k}{1-r_k}
,
\end{eqnarray*}
hence
$$
p_{\cdot 2}^B
=
\frac{r_1}{4}
+
p_{32}^B
=
\frac{r_1}{4}
\Bigg(
1
+
\sum_{k=2}^K
\frac{r_k}{1-r_k}
\Bigg)
.
$$
Of course, one then has~$p_{33}^B=(1-r_1)/4-p_{32}^B$ and $p_{\cdot 3}^B=(1/2)-p_{\cdot 2}^B$, but these values will not be needed for our purposes. Note that, in contrast with Model~A, \textcolor{black}{here $p_{\cdot 2}^B=P[Yg_1]$ may be different from~$P[Eg_1](=r_1/2)$. We argue that this is natural in the setup considered: while names are obviously picked by parents, it seems spurious to assume that the mechanism they adopt to choose names will ensure that~$P[Yg_1]=P[Eg_1]$. For instance, if it turns out that~$n_1$ is a very popular name, then one would expect that~$P[Yg_1]<P[Eg_1]$, since parents in two-daughters families are likely to pick this popular name for the first child. Another advantage of Model~B is that arbitrarily high popularity is allowed: $r_1$ may assume any value in~$(0,1)$ in this model, whereas the constraint~$P[Yg_1]=P[Eg_1]$ inherent to Model~A excludes that~$r_1$ exceeds~$1/2$. Allowing for arbitrarily high popularity of a name is natural and attractive from a mathematical point of view (if not from a practical point of view).}

\section{Results}
\label{sec:results}

In any model satisfying Assumptions~\ref{Assum1}--\ref{Assum2}, the probability that a family has a boy given that it has a girl named~$n_1$ is 
$$
P[Eb\cup Yb|Eg_1\cup Yg_1]
=
\frac{P[(Eb\cup Yb)\cap(Eg_1\cup Yg_1)]}{P[Eg_1\cup Yg_1]}
=
\frac{(r_1/4)+(r_1/4)}{(r_1/2)+p_{\cdot 2}}
;
$$
see Table~\ref{Table1}. In Model~A, we have~$p_{\cdot 2}=p_{\cdot 2}^A=r_1/2$, which yields 

\begin{Theor}[Model~A]
\label{TheorModelA}
The probability that a family has a boy given that it has a girl named~$n_1$ is 
$
P[Eb\cup Yb|Eg_1\cup Yg_1]=1/2$, irrespective of the value of~$r_1$.  
\end{Theor}

This is the usual result that, if a family has a girl named~$n_1$ (and if two sisters may not be given the same name), then both genders are equally likely for the other child. 
We now turn to Model~B, for which the situation is very much different. Up to renumbering the names~$n_2,\ldots,n_K$, we may of course assume that~$r_2\geq r_3\geq \ldots\geq r_K$ (note that this does not impose anything on~$r_1$). The model is thus indexed by~$\mathcal{R}_K=\{r=(r_1,r_2,\ldots,r_K)\in\ \textcolor{black}{(}0,1\textcolor{black}{)}\ :\sum_{k=1}^Kr_k=1, \,r_2\geq r_3\geq \ldots \geq r_K\}$. For any~$r\in \mathcal{R}_K$, the probability that a family has a boy given that it has a girl named~$n_1$ is then 
\begin{equation}
\label{probapretheor}
P[Eb\cup Yb|Eg_1\cup Yg_1]
=
\frac{(r_1/4)+(r_1/4)}{(r_1/2)+p_{\cdot 2}^B}
=
\frac{2}{3+\sum_{k=2}^K \frac{r_k}{1-r_k}}
\cdot
\end{equation}
As~$r\to (0,1,0,\ldots,0)$, this probability converges to zero. Since~$r_k>0$ for any~$k$, it trivially holds that this probability is strictly smaller than~$2/3$, an upper bound that is obtained as~$r\to (1,0,\ldots,0)$. From continuity, the probability in~(\ref{probapretheor}) can take any value in the interval~$\textcolor{black}{(}0,\frac{2}{3}\textcolor{black}{)}$. We therefore proved the following result.

\begin{Theor}[Model~B]
\label{TheorMainProba}
Fix an integer~$K\geq 2$. Then, for any~$r\in \mathcal{R}_K$, the probability that the family has a boy given that it has a girl named~$n_1$ is 
$$
p(r)
=
\frac{2}{3+\sum_{k=2}^K \frac{r_k}{1-r_k}}
\cdot
$$	
Moreover, $p(\mathcal{R}_K)$, the image of~$\mathcal{R}_K$ under the mapping~$p$, is~$\textcolor{black}{(}0,\frac{2}{3}\textcolor{black}{)}$.  
\end{Theor}

\begin{figure}[t]
\label{Fig1}
\begin{center}
	\includegraphics[width=110mm]{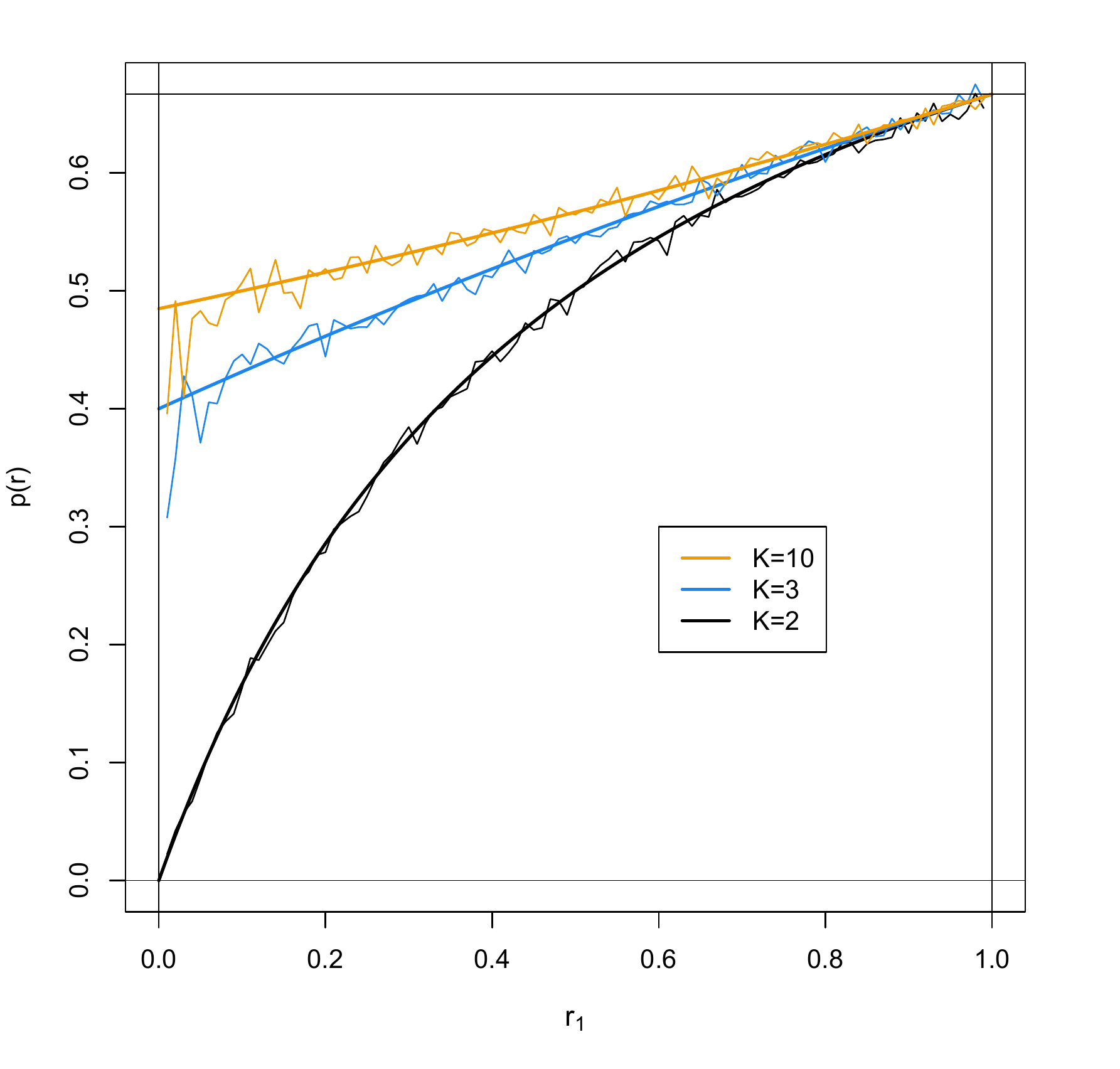}
\end{center}
\caption{Plots of~$p(r)$, with~$r=(r_1,(1-r_1)/(K-1),\ldots,(1-r_1)/(K-1))$, as a function of~$r_1$, for~$K=2,3$ and~$10$. The wiggly curves are obtained from Monte Carlo simulations, based on~$10,000$ replications, performed for any combination of~$r_1\in\{.01,.02,\ldots,.99\}$ and~$K\in\{2,3,10\}$.}
\end{figure}

For any~$K\geq 2$, the uniform configuration~$r=(1/K,\ldots,1/K)$ provides the case~$p(r)=1/2$ where both genders are equally likely for the other child. 
Figure~1
 plots~$p(r)$ as a function of~$r_1$ for various values of~$K$ in the framework where the names~$n_2,\ldots,n_K$ are equally likely (the figure also provides Monte Carlo simulation results that clearly support our expression of~$p(r)$ in Theorem~\ref{TheorMainProba}). It will be of interest to determine for which values of~$r_1$ there exists a corresponding configuration~$r=(r_1,\ldots,r_K)$ \textcolor{black}{for which the probability~$p(r)$ that the other child is a boy in the popularity Model~B agrees with the value, $1/2$, obtained in the traditional Model~A}. For this purpose, a crucial step is to  characterize, for each fixed~$K\geq 2$ and~$r_1\in \textcolor{black}{(}0,1\textcolor{black}{)}$, the values of~$p(r)$ that can be achieved. This is precisely what is done in the following result, whose proof is based on the Schur-concavity (for any~$r_1$) of the mapping~$(r_2,\ldots,r_K)\mapsto p(r_1,r_2,\ldots,r_K)$; see the appendix for a proof.  

\begin{Theor}[Model~B]
\label{TheorImagesr1}
Fix an integer~$K\geq 2$ and~$r_1\in\textcolor{black}{(}0,1\textcolor{black}{)}$. Define~$\mathcal{R}_K(c):=\{r\in\mathcal{R}_K:r_1=c\}$ the collection of values~$r$ for which~$r_1=c$. Then, for~$K=2$, 
$$
p(\mathcal{R}_K(r_1))
=
\bigg\{ 
\frac{2r_1}{2r_1+1}
\bigg\} 
,
$$
whereas, for~$K>2$, 
\begin{equation}
\label{Rr1K}
p(\mathcal{R}_K(r_1))
=
\textcolor{black}{\bigg(}
\frac{2r_1}{2r_1+1}
,
\frac{2(K-2+r_1)}{4(K-2+r_1)+1-Kr_1}
\bigg]
.
\end{equation}
For~$K>2$, the lower bound is obtained as~$r\to (r_1,1-r_1,0,\ldots,0)$, whereas the upper bound is achieved at~$r=(r_1,(1-r_1)/(K-1),\ldots,(1-r_1)/(K-1))$. 
\end{Theor}

The graphical illustration in 
Figure~2 
shows that, unless~$K=2$, typical popularity values for the name~$n_1$ (say, $r_1\leq 5\%$) will provide a rather wide range for the probability~$p(r)$ that the other child is a boy, that is, for such~$K$ and~$r_1$, this probability will much depend on the popularities of the remaining~$K-1$ names.
Note that the upper bound in~(\ref{Rr1K}) reduces to~$2r_1/(2r_1+1)$ for~$K=2$. As expected,  for any~$r_1\in\textcolor{black}{(}0,1\textcolor{black}{)}$, the (closures of the) feasible sets~$p(\mathcal{R}_K(r_1))$, $K=2,3,\ldots$, form a strictly increasing sequence with respect to inclusion, and
$$
\lim_{K\to\infty}
p(\mathcal{R}_K(r_1))
=
\textcolor{black}{\bigg(}
\frac{2r_1}{2r_1+1}
,
\frac{2}{4-r_1}
\bigg]
.
$$
Theorem~\ref{TheorImagesr1} also easily yields the following result.

\begin{Corol}[Model~B]
	\label{Coroll}
Fix an integer~$K\geq 2$ and define~$S_K:=\{r_1\in \textcolor{black}{(}0,1\textcolor{black}{)}: p(\mathcal{R}_K(r_1)) \ni \frac{1}{2}\}$, the collection of~$r_1$-values for which some configuration~$r=(r_1,\ldots,r_K)$ makes,  in a family having a girl named~$n_1$, both genders equally likely for the other child. Then, 
$$
S_2
=
\bigg\{
\frac{1}{2}
\bigg\}
\quad 
\textrm{ and }
\quad 
S_K
=
\bigg[
\frac{1}{K},\frac{1}{2}
\textcolor{black}{\bigg)}
\quad
\textrm{ for }K>2
.$$ 
\end{Corol}

\begin{figure}[t]
\label{Fig2}
\begin{center}
	\includegraphics[width=110mm]{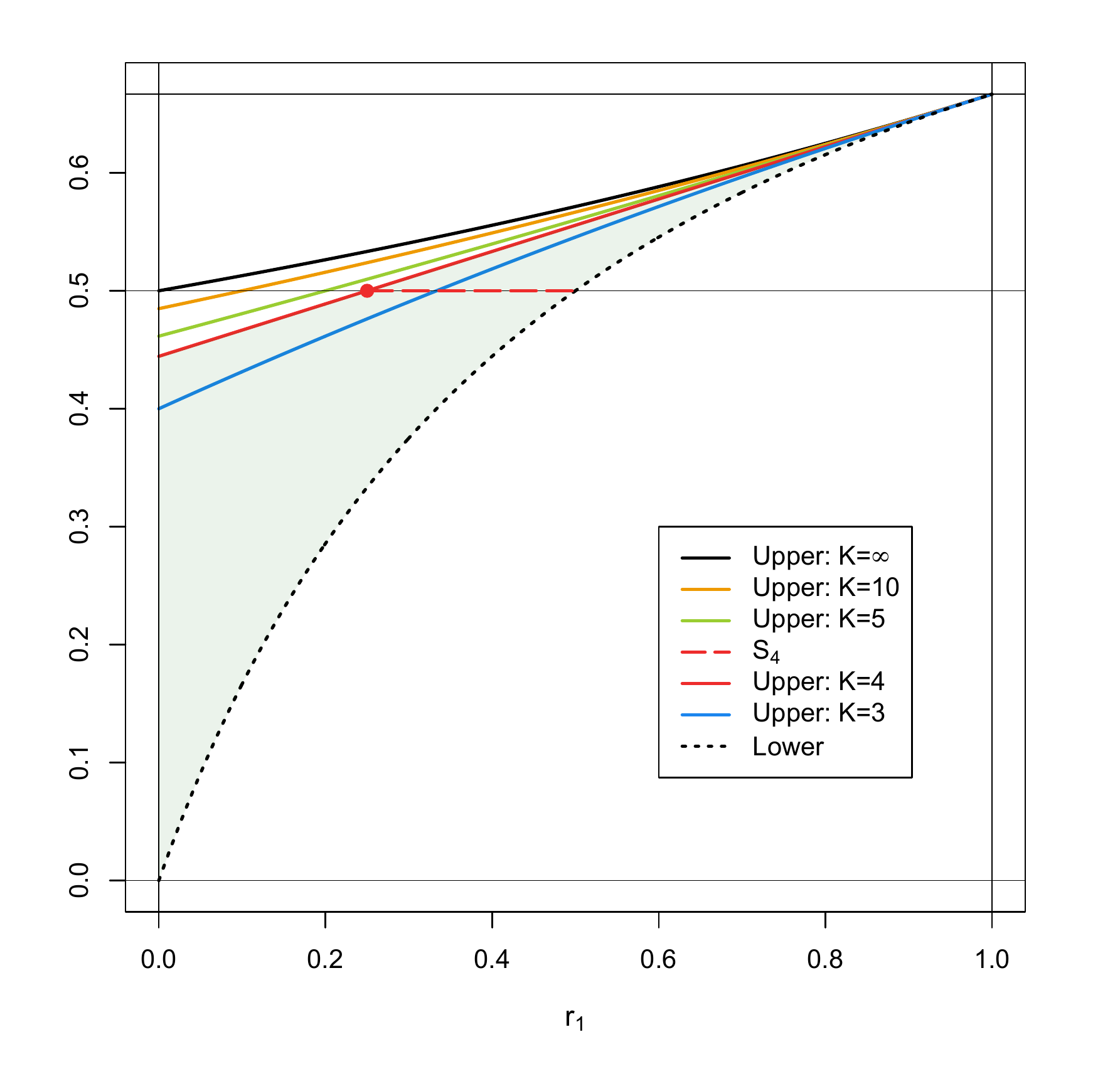}
\end{center}
\caption{Plots of the (fixed-$r_1$) upper bound of~$p(r)$ in~(\ref{Rr1K}) as a function of~$r_1$ for various values of~$K$, as well as the corresponding lower bound (that does not depend on~$K$). The shaded area thus emphasizes the feasible values of~$p(r)$ for~$K=4$. The dashed line segment corresponds to the set~$S_4$ in Corollary~\ref{Coroll}.}
\end{figure}

\begin{figure}[t]
\label{Fig3}
\begin{center}
	\includegraphics[width=110mm]{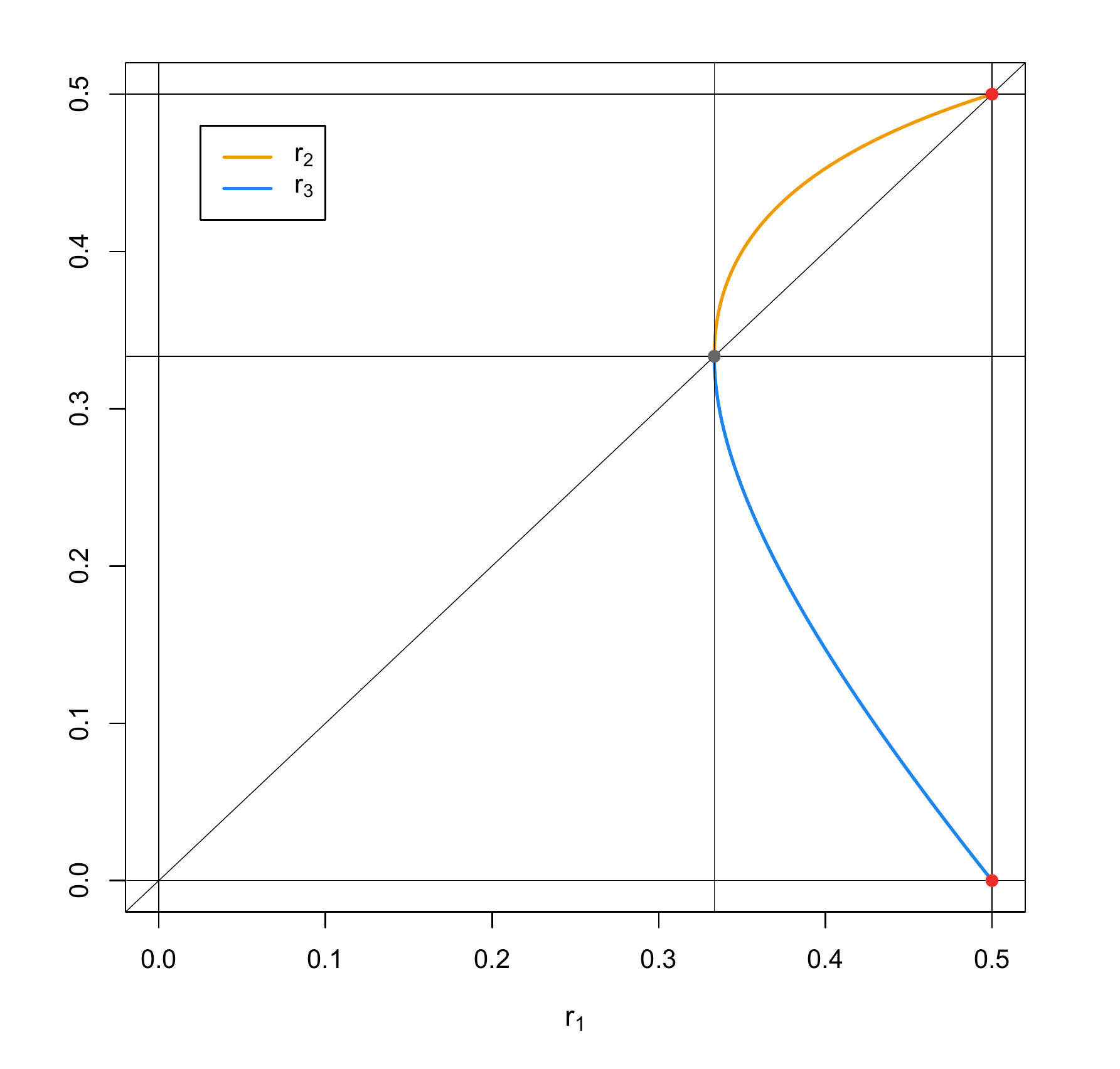}
\end{center}
\caption{This plots~$r_2$ and~$r_3$ as functions of~$r_1$ for all triples~$r=(r_1,r_2,r_3)\in\mathcal{R}_3$ providing~$p(r)=1/2$ in the case~$K=3$ (vertical and horizontal lines are plotted at~$1/3$ and~$1/2$). This confirms that~$S_3:=\{r_1\in\textcolor{black}{(}0,1\textcolor{black}{)}: p(\mathcal{R}_3(r_1)) \ni \frac{1}{2}\}$ is the interval~$[\frac{1}{3},\frac{1}{2}\textcolor{black}{)}$; see Corollary~\ref{Coroll}. It is seen that any solution provides~$r_1\in[r_2,r_3]$. }
\end{figure}

 As mentioned above, for any~$K\geq 2$, the uniform case~$r_1=r_2=\ldots=r_K(=1/K)$ provides~$p(r)=1/2$. For~$K=2$, this is the only case leading to equally likely genders for the second child. For~$K=3$, it is readily checked that the cases~$r=(r_1,r_2,r_3)\in\mathcal{R}_3$ providing~$p(r)=1/2$ are described by
$$
r_2
=
\frac{1-r_1}{2} 
+
\frac{\sqrt{(3r_1+1)^2-4}}{6} 
\quad
\textrm{and}
\quad
r_3 
= 
\frac{1-r_1}{2}
- 
\frac{\sqrt{(3r_1+1)^2-4}}{6} 
,
$$
with~$r_1\in [\frac{1}{3},\frac{1}{2}\textcolor{black}{)}$. 
Figure~\ref{Fig3}
 offers a graphical representation. Note that the results are in line with Corollary~\ref{Coroll}. \textcolor{black}{More generally, for~$K\geq 4$, the collections of~$(r_2,\ldots,r_K)$ making both genders equally likely for the second child is a manifold of dimension~$K-3$ in~$(\R^+_0)^{K-1}$.}

\section{Final comments}
\label{sec:finalcomments}

This paper revisits the name variant of the two-children problem and shows that, if it is known that the family has a girl named \emph{Florida}, say, then the probability that the other child is a boy may depend on the stochastic model that is adopted, even if one restricts to models that do not allow two sisters to have the same name. We show this by recalling that this probability is always~$1/2$ in the traditional model and by defining an alternative, natural, model in which this probability may assume any value in 
$$
\textcolor{black}{\bigg(}
\frac{2r_1}{2r_1+1}
,
\frac{2(K-2+r_1)}{4(K-2+r_1)+1-Kr_1}
\bigg]
,
$$
where~$r_1$ is the popularity of the name~\emph{Florida} (more precisely, $r_1$ is the probability that the first girl born in a family is given this name) and where~$K$ is the number of possible names (this assumes that~$K\geq 3$; for~$K=2$, this probability must be equal to~$2r_1/(2r_1+1)$). In this new model, the dependence of this probability on~$r_1$ clearly results from the heterogeneous way names are assigned to a first girl or a second one in a family.

It is of course natural to investigate whether or not the assumptions from Section~\ref{sec:models} are valid. For Assumption~\ref{Assum1}, this was discussed in \cite{Carlton2005}, where the authors conclude that neither Part~(i) nor Part~(ii) of this assumption actually holds in \textcolor{black}{practice}: more precisely, births of boys are more frequent than births of girls, and independence of genders is violated; in line with this, the null hypothesis that the number of boys follows a binomial distribution with parameters~$2$ and~$p$ for some unspecified~$p\in[0,1]$ is rejected at all usual significance levels (the p-value is below~$10^{-5}$). While results of the present note can be trivially extended to asymmetric gender probabilities, it is unclear how to deal with violation of the independence assumption (to the best of the authors' knowledge, this issue actually has not been touched for the classical Model~A). \textcolor{black}{While it would also} be of interest to test whether or not Assumption~\ref{Assum2} holds in practice, \textcolor{black}{it would be more urgent, in the context of the present paper, to focus on Assumption~3 and to investigate which model, among Model~A and Model~B, provides a better description of the real world. This could be done, on the basis of suitable data, by studying whether or not the proportion of girls named~$n_1$ among elder girls is different from the corresponding proportion among younger girls (a significant difference would make it necessary to consider Model~B; see  Section~\ref{sec:models}).}

\textcolor{black}{From} an inferential point of view, it is natural to estimate the popularity parameters~$r_k$, $k=1,\ldots,K$, by the observed frequencies of the various names among the collection of first girls born in two-children families (using the terminology adopted in the present note, first girls include elder girls as well as younger girls having an older brother). \textcolor{black}{This may not be the optimal solution, however, as there are likely ways to exploit} information among younger girls having an older sister, too. Yet, restricting to a subsample is a reasonable approach, and, as a matter of fact, it was also considered in \cite{Carlton2005}, where, after rejecting the null hypothesis that genders of both children are independent, the probability that a newborn is a boy is estimated by restricting to elder children.


\appendix

\section{Technical proofs}
\label{sec:appendix}

In this appendix, we prove Theorem~\ref{TheorImagesr1} and Corollary~\ref{Coroll}. 
\vspace{3mm}

\noindent
{\sc Proof of Theorem~\ref{TheorImagesr1}.}
The result for~$K=2$ trivially follows from Theorem~\ref{TheorMainProba}. Fix then~$K>2$ and~$r_1\in \textcolor{black}{(}0,1\textcolor{black}{)}$. Since the function~$x\mapsto g(x)=x/(1-x)$ is strictly convex over~$\textcolor{black}{(}0,1\textcolor{black}{)}$, Proposition~C.1a in Page~92 from \cite{Mar2011} yields that the function 
$$
(r_2,\ldots,r_K)\mapsto 
\sum_{k=2}^K \frac{r_k}{1-r_k}
$$	
is strictly Schur-convex on~$\textcolor{black}{(}0,1\textcolor{black}{)}^{K-1}$, so that
$$
(r_2,\ldots,r_K)
\mapsto 
p(r_1,r_2,\ldots,r_K)
=
\frac{2}{3+\sum_{k=2}^K \frac{r_k}{1-r_k}}
$$	
is strictly Schur-concave on~$\textcolor{black}{(}0,1\textcolor{black}{)}^{K-1}$. For any~$r\in \mathcal{R}_K(r_1)$, we have~$\sum_{k=2}^K r_k=1-r_1$, hence also
$$
\bigg(\frac{1-r_1}{K-1},\ldots,\frac{1-r_1}{K-1}\Bigg)
\prec
(r_2,\ldots,r_K)
\prec
(1-r_1,0,\ldots,0)
,
$$
where~$\prec$ refers to the usual majorization ordering; see, e.g., pages~8--9 in \cite{Mar2011}. The above Schur-concavity therefore implies that
$$  
p(r_1,1-r_1,0,\ldots,0)
\leq
p(r)
\leq
p\bigg(r_1,\frac{1-r_1}{K-1},\ldots,\frac{1-r_1}{K-1}\Bigg)
$$
for any~$r\in \mathcal{R}_K(r_1)$. Since direct evaluation provides
$$
p(r_1,1-r_1,0,\ldots,0)
=
\frac{2r_1}{2r_1+1}
$$
and
$$
p\bigg(r_1,\frac{1-r_1}{K-1},\ldots,\frac{1-r_1}{K-1}\Bigg)
=
\frac{2(K-2+r_1)}{4(K-2+r_1)+1-Kr_1}
,
$$
the result follows from continuity of~$r\mapsto p(r)$.
\cqfd
\vspace{3mm}


\noindent
{\sc Proof of Corollary~\ref{Coroll}.}
The result trivially holds for~$K=2$, so we may restrict to~$K>2$. From Theorem~\ref{TheorImagesr1}, we have that~$r_1\in S_K$ --- i.e., $\frac{1}{2}\in p(\mathcal{R}_K(r_1))$ --- if and only if
$$
\frac{2r_1}{2r_1+1}
<
\frac{1}{2}
\quad
\textrm{ and }
\quad
\frac{1}{2}
\leq
\frac{2(K-2+r_1)}{4(K-2+r_1)+1-Kr_1}
,
$$
that is, if and only if
$$
r_1
<
\frac{1}{2}
\quad
\textrm{ and }
\quad
r_1
\geq 
\frac{1}{K}
,
$$
which establishes the result. 
\cqfd
\vspace{3mm}

\vspace{7mm}

 \noindent \textbf{Acknowledgements}
\vspace{2mm}
 
\noindent
This research is supported by the Program of Concerted Research Actions (ARC) of the Universit\'{e} libre de Bruxelles. This note results from exchanges following a talk of the Alta\"{i}r conference cycle in Brussels; the authors would like to thank the organisers. 
\vspace{0mm}

\bibliographystyle{chicago}
\bibliography{Paper.bib}

\end{document}